\documentclass[12pt]{amsart}
\usepackage{fullpage}
\usepackage{url}

\newtheorem{theorem}{Theorem}

\newtheorem{lemma}[theorem]{Lemma}

\theoremstyle{remark}

\begin{document}

\title{Further applications of a power series method for pattern avoidance}

\author{Narad Rampersad}

\address{Department of Mathematics and Statistics \\
University of Winnipeg \\
515 Portage Avenue \\
Winnipeg, Manitoba R3B 2E9 (Canada)}

\email{n.rampersad@uwinnipeg.ca}

\thanks{The author is supported by an NSERC Postdoctoral
Fellowship.}


\date{\today}

\begin{abstract}
In combinatorics on words, a word $w$ over an alphabet $\Sigma$ is
said to avoid a pattern $p$ over an alphabet $\Delta$ if there is no
factor $x$ of $w$ and no non-erasing morphism $h$ from $\Delta^*$ to
$\Sigma^*$ such that $h(p) = x$.  Bell and Goh have recently applied
an algebraic technique due to Golod to show that for a certain wide
class of patterns $p$ there are exponentially many words of length $n$
over a $4$-letter alphabet that avoid $p$.  We consider some further
consequences of their work.  In particular, we show that any pattern
with $k$ variables of length at least $4^k$ is avoidable on the binary
alphabet.  This improves an earlier bound due to Cassaigne and Roth.
\end{abstract}

\maketitle

\section{Introduction}

In combinatorics on words, the notion of an avoidable/unavoidable
pattern was first introduced (independently) by Bean, Ehrenfeucht, and
McNulty \cite{BEM79} and Zimin \cite{Zim84}.  Let $\Sigma$ and
$\Delta$ be alphabets: the alphabet $\Delta$ is the \emph{pattern
  alphabet} and its elements are \emph{variables}.  A \emph{pattern}
$p$ is a non-empty word over $\Delta$.  A word $w$ over $\Sigma$ is an
\emph{instance of $p$} if there exists a non-erasing morphism
$h:\Delta^*\to\Sigma^*$ such that $h(p) = w$.  A pattern $p$ is
\emph{avoidable} if there exists infinitely many words $x$ over a
finite alphabet such that no factor of $x$ is an instance of $p$.
Otherwise, $p$ is \emph{unavoidable}.  If $p$ is avoided by infinitely
many words on an $m$-letter alphabet then it is said to be
\emph{$m$-avoidable}.  The survey chapter in Lothaire
\cite[Chapter~3]{Lot02} gives a good overview of the main results
concerning avoidable patterns.

The classical results of Thue \cite{Thu06,Thu12} established that the
pattern $xx$ is $3$-avoidable and the pattern $xxx$ is $2$-avoidable.
Schmidt \cite{Sch89} (see also \cite{Ram08}) proved that any binary
pattern of length at least $13$ is $2$-avoidable; Roth \cite{Rot92}
showed that the bound of $13$ can be replaced by $6$.  Cassaigne
\cite{Cas93} and Vani\u{c}ek \cite{Van89} (see \cite{GV91})  determed
exactly the set of binary patterns that are $2$-avoidable.

Bean, Ehrenfeucht, and McNulty \cite{BEM79} and Zimin \cite{Zim84}
characterized the avoidable patterns in general.  Let us call a
pattern $p$ for which all variables occurring in $p$ occur at least
twice a \emph{doubled pattern}.  A consequence of the characterization
of the avoidable patterns is that any doubled pattern is avoidable.
Bell and Goh \cite{BG07} proved the much stronger result that every
doubled pattern is $4$-avoidable.  Cassaigne and Roth (see
\cite{Cas94} or \cite[Chapter~3]{Lot02}) proved that any pattern
containing $k$ distinct variables and having length greater than
$200\cdot 5^k$ is $2$-avoidable.  In this note we apply the arguments
of Bell and Goh to show the following result, which improves that of
Cassaigne and Roth.

\begin{theorem}\label{main}
Let $k$ be a positive integer and let $p$ be a pattern containing
$k$ distinct variables.
\begin{itemize}
\item[(a)] If $p$ has length at least $2^k$ then $p$ is $4$-avoidable.
\item[(b)] If $p$ has length at least $3^k$ then $p$ is $3$-avoidable.
\item[(c)] If $p$ has length at least $4^k$ then $p$ is $2$-avoidable.
\end{itemize}
\end{theorem}

\section{A power series approach}

Rather than simply wishing to show the avoidability of a
pattern $p$, one may wish instead to determine the number of
words of length $n$ over an $m$-letter alphabet that avoid $p$
(see, for instance, Berstel's survey \cite{Ber05}).  Brinkhuis
\cite{Bri83} and Brandenburg \cite{Bra83} showed that there are
exponentially many words of length $n$ over a $3$-letter alphabet that
avoid the pattern $xx$.  Similarly, Brandenburg showed that there are
exponentially many words of length $n$ over a $2$-letter alphabet that
avoid the pattern $xxx$. 

As previously mentioned, Bell and Goh proved that every doubled
pattern is $4$-avoidable.  In fact, they proved the stronger result
that there are exponentially many words of length $n$ over a
$4$-letter alphabet that avoid a given doubled pattern.  Their main
tool in obtaining this result is the following.

\begin{theorem}[Golod]\label{golod}
Let $S$ be a set of words over an $m$-letter alphabet, each word of
length at least $2$.  Suppose that for each $i \geq 2$, the set $S$
contains at most $c_i$ words of length $i$.  If the power series
expansion of
\begin{equation}\label{pow_ser}
G(x) := \left( 1 - mx + \sum_{i \geq 2} c_i x^i \right)^{-1}
\end{equation}
has non-negative coefficients, then there are least $[x^n]G(x)$
words of length $n$ over an $m$-letter alphabet that avoid $S$.
\end{theorem}

Theorem~\ref{golod} was originally presented by Golod (see
Rowen \cite[Lemma~6.2.7]{Row88}) in an algebraic setting.  We
have restated it here using combinatorial terminology.  The proof
given in Rowen's book also is phrased in algebraic terminology; in
order to make the technique perhaps a little more accessible to
combinatorialists, we present a proof of Theorem~\ref{golod} using
combinatorial language.

\begin{proof}[Proof of Theorem~\ref{golod}]
For two power series $f(x) = \sum_{i \geq 0} a_i x^i$ and $g(x) =
\sum_{i \geq 0} b_i x^i$, we write $f \geq g$ to mean that $a_i \geq
b_i$ for all $i \geq 0$.  Let $F(x) := \sum_{i \geq 0} a_i x^i$, where
$a_i$ is the number of words of length $i$ over an $m$-letter alphabet
that avoid $S$.  Let $G(x) = \sum_{i \geq 0} b_i x^i$ be the power series
expansion of $G$ defined above.  We wish to show $F \geq G$.

For $k \geq 1$, there are $m^k - a_k$ words $w$ of length $k$ over an
$m$-letter alphabet that contain a word in $S$ as a factor.  On the
other hand, for any such $w$ either (a) $w = w'a$, where $a$ is a
single letter and $w'$ is a word of length $k-1$ containing a word in
$S$ as a factor; or (b) $w = xy$, where $x$ is a word of length $k-j$
that avoids $S$ and $y \in S$ is a word of length $j$.  There are at
most $(m^{k-1} - a_{k-1})m$ words $w$ of the form (a), and there are
at most $\sum_j a_{k-j}c_j$ words $w$ of the form (b).  We thus have
the inequality
\[
m^k - a_k \leq (m^{k-1} - a_{k-1})m + \sum_j a_{k-j}c_j.
\]
Rearranging, we have
\begin{equation}\label{ineq}
a_k - a_{k-1}m + \sum_j a_{k-j}c_j \geq 0,
\end{equation}
for $k \geq 1$.

Consider the function
\begin{eqnarray*}
H(x) & := & F(x) \left(1 - mx + \sum_{j \geq 2} c_j x^j\right) \\
& = & \left(\sum_{i \geq 0} a_i x^i\right)\left(1 - mx + \sum_{j \geq 2} c_j x^j\right).
\end{eqnarray*}
Observe that for $k \geq 1$, we have $[x^k]H(x) = a_k - a_{k-1}m +
\sum_j a_{k-j}c_j$.  By (\ref{ineq}), we have $[x^k]H(x) \geq 0$ for
$k \geq 1$.  Since $[x^0]H(x) = 1$, the inequality $H \geq 1$ holds,
and in particular, $H-1$ has non-negative coefficients.  We conclude
that $F = HG = (H-1)G + G \geq G$, as required.
\end{proof}

Theorem~\ref{golod} bears a certain resemblance to the
Goulden--Jackson cluster method \cite[Section~2.8]{GJ04}, which also
produces a formula similar to (\ref{pow_ser}).  The cluster method
yields an exact enumeration of the words avoiding the set $S$ but
requires $S$ to be finite.  By constrast, Theorem~\ref{golod} only
gives a lower bound on the number of words avoiding $S$, but now the
set $S$ can be infinite.

Theorem~\ref{golod} can be viewed as a non-constructive method to show
the avoidability of patterns over an alphabet of a certain size.  In
this sense it is somewhat reminiscent of the probabilistic approach to
pattern avoidance using the Lov\'asz local lemma (see \cite{Bec81,
  Cur05}).  For pattern avoidance it may even be more powerful than
the local lemma in certain respects.  For instance, Pegden
\cite{Peg09} proved that doubled patterns are $22$-avoidable using the
local lemma, whereas Bell and Goh were able to show $4$-avoidability
using Theorem~\ref{golod}.  Similarly, the reader may find it a
pleasant exercise to show using Theorem~\ref{golod} that there are
infinitely many words avoiding $xx$ over a $7$-letter alphabet; as far
as we are aware, the smallest alphabet size for which the avoidability
of $xx$ has been shown using the local lemma is $13$ \cite{Sha08}.

\section{Proof of Theorem~\ref{main}}

To prove Theorem~\ref{main} we begin with some lemmas.

\begin{lemma}\label{multiple}
Let $k \geq 1$ and $m \geq 2$ be integers.  If $w$ is a word of
length at least $m^k$ over a $k$-letter alphabet, then $w$ contains
a non-empty factor $w'$ such that the number of occurrences of each
letter in $w'$ is a multiple of $m$.
\end{lemma}

\begin{proof}
Suppose $w$ is over the alphabet $\Sigma = \{1,2,\ldots,k\}$.  Define
the map $\psi : \Sigma^* \to \mathbb{N}^k$ that maps a word $x$ to the
$k$-tuple $[|x|_1 \bmod m,\ldots,|x|_k \bmod m]$, where $|x|_a$
denotes the number of occurrences of the letter $a$ in $x$.  For each
prefix $w_i$ of length $i$ of $w$, let $v_i = \psi(w_i)$.  Since $w$
has length at least $m^k$, $w$ has at least $m^k + 1$ prefixes, but
there are at most $m^k$ distinct tuples $v_i$.  There exists therefore
$i<j$ such that $v_i=v_j$.  However, if $w'$ is the suffix of $w_j$ of
length $j-i$, then $\psi(w') = v_j - v_i = [0,\ldots,0]$, and hence
the number of occurrences of each letter in $w'$ is a multiple of $m$.
\end{proof}

\begin{lemma}[\cite{BG07}]\label{match}
Let $k \geq 1$ be a integer and let $p$ be a pattern over
the pattern alphabet $\{x_1, \ldots, x_k\}$.  Suppose that
for $1 \leq i \leq k$, the variable $x_i$ occurs $a_i \geq 1$ times
in $p$.  Let $m \geq 2$ be an integer and let $\Sigma$ be an
$m$-letter alphabet.  Then for $n\geq 1$, the number of words of
length $n$ over $\Sigma$ that are instances of the pattern $p$ is at most
$[x^n]C(x)$, where
\[
 C(x) := \sum_{i_1 \geq 1} \cdots \sum_{i_k
\geq 1} m^{i_1+\cdots+i_k} x^{a_1i_1+\cdots+a_ki_k}.
\]
\end{lemma}

For the proof of the next lemma, we essentially follow the approach of Bell
and Goh.

\begin{lemma}\label{growth}
Let $k \geq 2$ be an integer and let $p$ be a pattern over a
$k$-letter pattern alphabet such that every variable occuring in $p$
occurs at least $\mu$ times.
\begin{itemize}
\item[(a)] If $\mu=3$, then for $n \geq 0$, there are at least
$2.94^n$ words of length $n$ avoiding $p$ over a $3$-letter alphabet.
\item[(b)] If $\mu=4$, then for $n \geq 0$, there are at least
$1.94^n$ words of length $n$ avoiding $p$ over a $2$-letter alphabet.
\end{itemize}
\end{lemma}

\begin{proof}
Let $(m,\mu) \in \{(3,3),(2,4)\}$ and let $\Sigma$ be an $m$-letter
alphabet.  Define $S$ to be the set of all words over $\Sigma$ that
are instances of the pattern $p$.  By Lemma~\ref{match}, the number
of words of length $n$ in $S$ is at most $[x^n]C(x)$, where
\[
C(x) := \sum_{i_1 \geq 1} \cdots \sum_{i_k
\geq 1} m^{i_1+\cdots+i_k} x^{a_1i_1+\cdots+a_ki_k},
\]
and for $1 \leq i \leq k$ we have $a_i \geq \mu$.  Define
\[
B(x) := \sum_{i \geq 0} b_ix^i = (1 - mx +C(x))^{-1},
\]
and set $\lambda := m-0.06$. We claim that $b_n \geq \lambda b_{n-1}$
for all $n \geq 0$.  This suffices to prove the lemma, as we would
then have $b_n \geq \lambda^n$ and the result follows by an
application of Theorem~\ref{golod}.

We prove the claim by induction on $n$.  When $n=0$, we have $b_0=1$ and
$b_1=m$.  Since $m > \lambda$, the inequality $b_1 \geq \lambda b_0$
holds, as required.  Suppose that for all $j<n$, we have $b_j \geq
\lambda b_{j-1}$.  Since $B = (1 - mx +C)^{-1}$, we have
$B(1 - mx +C) = 1$.  Hence $[x^n]B(1-mx+C) = 0$ for $n \geq 1$.
However,
\[
B(1-mx+C) = \left(\sum_{i \geq 0} b_ix^i\right) \left(1 - mx +
\sum_{i_1 \geq 1} \cdots \sum_{i_k \geq 1} m^{i_1+\cdots+i_k}
x^{a_1i_1+\cdots+a_ki_k}\right),
\]
so
\[
[x^n]B(1-mx+C) = b_n - b_{n-1}m + \sum_{i_1 \geq 1} \cdots \sum_{i_k
\geq 1} m^{i_1+\cdots+i_k} b_{n-(a_1i_1+\cdots+a_ki_k)} = 0.
\]
Rearranging, we obtain
\[
b_n = \lambda b_{n-1} + (m-\lambda)b_{n-1} - \sum_{i_1 \geq 1} \cdots \sum_{i_k
\geq 1} m^{i_1+\cdots+i_k} b_{n-(a_1i_1+\cdots+a_ki_k)}.
\]
To show $b_n \geq \lambda b_{n-1}$ it therefore suffices to show
\begin{equation}\label{bn_growth}
(m-\lambda)b_{n-1} - \sum_{i_1 \geq 1} \cdots \sum_{i_k
\geq 1} m^{i_1+\cdots+i_k} b_{n-(a_1i_1+\cdots+a_ki_k)} \geq 0.
\end{equation}
Since $b_j \geq \lambda b_{j-1}$ for all $j<n$, we have $b_{n-i} \leq
b_{n-1}/\lambda^{i-1}$ for $1 \leq i \leq n$.  Hence
\begin{eqnarray*}
& & \sum_{i_1 \geq 1} \cdots \sum_{i_k
\geq 1} m^{i_1+\cdots+i_k} b_{n-(a_1i_1+\cdots+a_ki_k)} \\
& \leq & \sum_{i_1 \geq 1} \cdots \sum_{i_k
\geq 1} m^{i_1+\cdots+i_k} \frac{\lambda
b_{n-1}}{\lambda^{a_1i_1+\cdots+a_ki_k}} \\
& = & \lambda b_{n-1} \sum_{i_1 \geq 1} \cdots \sum_{i_k
\geq 1} \frac{m^{i_1+\cdots+i_k}}{\lambda^{a_1i_1+\cdots+a_ki_k}} \\
& = & \lambda b_{n-1} \sum_{i_1 \geq 1}
\frac{m^{i_1}}{\lambda^{a_1i_1}} \cdots \sum_{i_k
\geq 1} \frac{m^{i_k}}{\lambda^{a_ki_k}} \\
& \leq & \lambda b_{n-1} \sum_{i_1 \geq 1}
\frac{m^{i_1}}{\lambda^{\mu i_1}} \cdots \sum_{i_k
\geq 1} \frac{m^{i_k}}{\lambda^{\mu i_k}} \\
& = & \lambda b_{n-1} \left( \sum_{i \geq 1} \frac{m^i}{\lambda^{\mu
      i}} \right)^k \\
& = & \lambda b_{n-1} \left( \frac{m/\lambda^\mu}{1-m/\lambda^\mu}
\right)^k \\
& = & \lambda b_{n-1} \left( \frac{m}{\lambda^\mu - m} \right)^k \\
& \leq & \lambda b_{n-1} \left( \frac{m}{\lambda^\mu - m} \right)^2.
\end{eqnarray*}

In order to show that (\ref{bn_growth}) holds, it thus suffices to
show that
\[
m-\lambda \geq \lambda \left( \frac{m}{\lambda^\mu - m} \right)^2.
\]
Recall that $m-\lambda = 0.06$.  For $(m,\mu) = (3,3)$ we have
\[
2.94 \left( \frac{3}{2.94^3 - 3} \right)^2 = 0.052677\cdots \leq 0.06,
\]
and for $(m,\mu) = (2,4)$ we have
\[
1.94 \left( \frac{2}{1.94^4 - 2} \right)^2 = 0.052439\cdots \leq 0.06,
\]
as required.  This completes the proof of the inductive claim and the
proof of the lemma.
\end{proof}

We can now complete the proof of Theorem~\ref{main}.  Let $p$ be a
pattern with $k$ variables.  If $p$ has length at least $2^k$, then by
Lemma~\ref{multiple}, the pattern $p$ contains a non-empty factor $p'$
such that each variable occurring in $p'$ occurs at least twice.
However, Bell and Goh showed that such a $p'$ is $4$-avoidable and
hence $p$ is $4$-avoidable.

Similarly, if $p$ has length at least $3^k$ (resp. $4^k)$, then by
Lemma~\ref{multiple}, the pattern $p$ contains a non-empty factor $p'$
such that each variable occurring in $p'$ occurs at least $3$ times
(resp. $4$ times).  If $p'$ contains only one distinct variable,
recall that we have already noted in the introduction that the pattern
$xxx$ is $2$-avoidable (and hence also $3$-avoidable).  If $p'$
contains at least two distinct variables, then by Lemma~\ref{growth},
the pattern $p'$ is $3$-avoidable (resp. $2$-avoidable), and hence the
pattern $p$ is $3$-avoidable (resp. $2$-avoidable).  This completes
the proof of Theorem~\ref{main}.

Recall that Cassaigne and Roth showed that any pattern $p$ over $k$
variables of length greater than $200\cdot 5^k$ is $2$-avoidable.
Their proof is constructive but is rather difficult.  We are able to
obtain the much better bound of $4^k$ non-constructively by a somewhat
simpler argument.  Cassaigne suggests (see the open problem
\cite[Problem~3.3.2]{Lot02}) that the bound of $3^k$ in
Theorem~\ref{main}(b) can perhaps be replaced by $2^k$ and that the
bound of $4^k$ in Theorem~\ref{main}(c) can perhaps be replaced by
$3\cdot 2^k$.  Note that the bound of $2^k$ in Theorem~\ref{main}(a)
is optimal, since the Zimin pattern on $k$-variables (see
\cite[Chapter~3]{Lot02}) has length $2^k-1$ and is unavoidable.

\section*{Acknowledgments}

We thank Terry Visentin for some helpful discussions concerning
Theorem~\ref{golod} and the Goulden--Jackson cluster method.

\end{document}